\documentclass[12pt,reqno]{amsart}
\usepackage{amsthm,amsfonts,amssymb,euscript,mathrsfs}

\def\a{{\alpha}}
\def\b{{\beta}}

\def\d{\delta}

\def\e{\epsilon}

\def\hat{\widehat}

\def\FF{{\mathcal F}}

\def\R{{\mathbb{R}}}
\def\C{{\mathbb{C}}}

\def\S{{Schr\"{o}dinger }}

\def\Z{{\mathbb{Z}}}

\def\l{\langle}
\def\r{\rangle}

\newtheorem{theorem}{Theorem}[section]
\newtheorem{lemma}[theorem]{Lemma}
\newtheorem{proposition}[theorem]{Proposition}

\newtheorem{remark}[theorem]{Remark}

\newtheorem{convention}{Convention}

\numberwithin{equation}{section}

\begin{document}

\title[Scattering for the Zakharov system]{Scattering for the Zakharov system in 3 dimensions}

\author{Zaher Hani, \;Fabio Pusateri \;\& \;Jalal Shatah}


\begin{abstract}
We prove global existence and scattering for small localized solutions of the Cauchy problem for the Zakharov system in $3$ space dimensions.
The wave component is shown to decay pointwise at the optimal rate of $t^{-1}$, 
whereas the \S component decays almost at a rate of $t^{-7/6}$.
\end{abstract}

\maketitle

\section{Introduction}
The purpose of this manuscript is to study the asymptotic behavior of small solutions of the Zakharov system given by:
\begin{equation}
\label{Z}
\tag{Z}
\left\{
\begin{array}{l}
i\partial_t u + \Delta u =  n u
\\
\\
\Box n = \Delta |u|^2 \, ,
\end{array}
\right.
\end{equation}
where $(u,n) : (t,x) \in \R \times \R^3 \rightarrow \C \times \R$, and the initial data are taken to be:
\begin{equation*}
 u(0,x) = u_0(x) \, , \quad n(0,x) = n_0(x) \, , \quad \partial_t n(0,x) = n_1(x) \,.
\end{equation*}
Our main result is the following:
\begin{theorem}
\label{maintheo} 
Suppose that the initial data  $u_0$, $n_0$, $n_1$ satisfy:
\begin{align}
\label{data}
& {\| u_0 \|}_{H^{N+1}} + {\| {\langle x \rangle}^2 u_0 \|}_{L^2} \leq \e_0
\\
&\label{data2} {\left\| \left( \Lambda n_0, n_1 \right) \right\|}_{H^{N-1}} 
+ {\left\| \langle \Lambda \rangle \left( \Lambda n_0,  n_1 \right) \right\|}_{\dot{B}^0_{1,1}}
+ {\left\| \langle x \rangle \left( n_0, {\langle x \rangle} n_1 \right) \right\|}_{H^1} \leq \e_0 \, ,
\end{align}
for some small $\epsilon_0$ and some large integer $N$. Then  the Cauchy problem for the Zakharov system \eqref{Z} 
admits a unique global solution such that
\begin{equation}
{\| u(t) \|}_{L^\infty} \lesssim \frac{\epsilon_0}{t^{7/6 -}} \quad , \quad {\| n(t) \|}_{L^\infty} \lesssim \frac{\epsilon_0}{t} \, .
\end{equation}
As a consequence the solution $( u(t), n(t) )$ scatters to a linear solution as $t \rightarrow \infty$.
\end{theorem}
Here $\langle x \rangle$ is used to denote $\sqrt{1+|x|^2}$, 
$\Lambda:=\sqrt{-\Delta}$, and the definition of the Besov norm $\dot B^0_{1,1}$ is recalled in \eqref{Besov}.

The system \eqref{Z} is an important model in plasma physics and has been under intensive investigation by physicists and mathematicians.
It was derived by V. Zakharov in \cite{Zakharov} to model Langmuir
waves in plasma, in which context it describes the interaction between a high-frequency electromagnetic wave (the \S component $u$) 
with an acoustic wave (the component $n$) \cite{HamPlasma}. 
It serves as a simplified model for laser-plasma interaction where the function $u$ stands for the complex envelope of the electric field 
whereas $n$ stands for the mean density fluctuation of ions or electrons from the unperturbed plasma density. 
The Schr\"odinger operator appears as a three-scale approximation of 
Maxwell's equations whereas the wave component is a long-wave approximation of the Euler equations 
in the fundamental Euler-Maxwell system governing plasma motion. 
We refer to \cite{HamPlasma, SulemBook,Texier} and references therein for more background on the physical significance of \eqref{Z}.

From the mathematical side, there has been considerable work on local and global well-posedness of solutions with rough data 
through the works of Kenig, Ponce and Vega \cite{KPV1}, Bourgain and Colliander \cite{BC},  Ginibre, Tsutsumi and Velo \cite{GTV},
Bejenaru, Herr, Holmer and Tataru \cite{BHHT} and Bejenaru and Herr \cite{BH} 
(cf. the references in the cited works 
for previous well-posedness results).
Global well-posedness for small data in the energy space is obtained in \cite{BC} by combining local well-posedness and conservation laws. 
Many works have also dealt with singular limits related to the Zakharov system
and with the rigorous derivation of the system in various limiting regimes from other equations and vice versa. 
We refer the reader to the work of Texier \cite{Texier}, where \eqref{Z} is derived from the Euler-Maxwell equations, 
and the work of Masmoudi and Nakanishi \cite{MN}, where it is obtained from the Klein-Gordon-Zakharov system
(cf. references in \cite{Texier,MN} for previous results in this vein).

Concerning the scattering question, most of the previous work has been carried out for the  final value problem, i.e. data at $t = \infty$,  
instead of the Cauchy problem, as in the papers of Ozawa and Tsutsumi  \cite{OT}, Shimomura \cite{ShimoZak} and Ginibre and Velo \cite{GV-Zak}.
Similar work on the final value problem has also been dedicated to other coupled systems of \S and wave equations,
like in the papers of Ginibre and Velo \cite{GV-WS1,GV-WS2,GV-WS3} and Shimomura \cite{ShimoWS,ShimoMS}.
It is interesting to notice that the Wave-\S system considered in \cite{GV-WS1,GV-WS2,GV-WS3} and \cite{ShimoWS}
has a similar structure to the Zakharov system \eqref{Z}: the equation for $u$ 
is identical while the nonlinearity in the equation for $n$ is ${|u|}^2$ instead of $\Delta{|u|}^2$.
In the case of the Wave-\S system {\it modified} wave operators for the final value problem are constructed in the previously cited papers.
In the case of the Zakharov system \eqref{Z} that we investigate here, we are able to show 
(linear) scattering for solutions of the Cauchy problem.
This is made possible because the nonlinearity $\Delta{|u|}^2$ 
possesses a null structure which gives stronger control on the solution.

The only work that deals with small-data scattering for the Cauchy problem of the Zakharov system (or any other Wave-\S system in $3$ dimensions) 
is an important recent manuscript by Guo and Nakanishi \cite{GN} where the authors consider small radial solutions in the energy space.
The assumption of radial symmetry leads to a wider range of Strichartz estimates 
(with better time integrability) which allows one to close an iteration argument for data in the energy space. 
The great advantage of working at this level of regularity, that is controlled\footnote{
At least in an appropriately defined coercive regime (cf. \cite{GNW}).} by conserved quantities,
is that it allows the authors to tackle the large-data problem (see Guo, Nakanishi and Wang \cite{GNW}). 
Nonetheless, it is important to emphasize that the proof in \cite{GN} relies heavily on the radial assumption
which allows access to linear estimates that do not hold in the non-radial setting.
As we argue below, the main difficulties in dealing with the scattering problem without the radial assumption 
are the slow decay of the wave component and the mixed linear parts, 
which simultaneously rule out the possibility of using vectorfields and create many resonances.

%
%
%
%
%

\section{Preliminary setup and outline of the proof}
To prove Theorem \ref{maintheo} we utilize the space time resonance method \cite{GMS1,GMS2} as our general framework, 
rely on some of the ideas used in \cite{PS}, 
and on decomposing the nonlinearities according to the real space support of the interacting profiles.
Writing $w_{\pm}=\Lambda^{-1} (i\partial_t\pm \Lambda) n $, then  system   \eqref{Z} becomes
\begin{equation}
\label{Z_0}
\tag{Z$_0$}
\left\{
\begin{array}{l}
i\partial_t u +\Delta u = \frac{1}{2}(w_+ u - w_- u)
\\
\\
i\partial_t w_\pm \mp \Lambda w_\pm = \Lambda {|u|}^2  \, .
\end{array}
\right.
\end{equation}
Let  $f = e^{-it\Delta} u$ and $g_\pm=e^{\pm it\Lambda} w_\pm$ denote the profiles, and let $\hat f = \FF f$ and $\hat g_\pm =\FF g_\pm$ denote  their Fourier transforms, then we have from Duhamel's formula\begin{subequations}
\label{eq:profile}
\begin{align}
\label{inteqf}
&\hat f(t,\xi) = \hat f(0,\xi) \mp \sum_{\pm} i\int_0^t \int_{\R^3} e^{is \phi_{\pm}(\xi,\eta)}\hat f(\xi-\eta,s) \hat g_\pm(\eta,s) d\eta ds
\\
\label{inteqg}
&\hat g_\pm(t,\xi)  = \hat g_\pm (0,\xi)  - i\int_0^t \int_{\R^3}|\xi| e^{is \psi_{\pm}(\xi,\eta)}
			\hat f(\xi-\eta,s) \overline{\hat{f}} (\eta,s) d\eta ds \, ,
\end{align}
\end{subequations}
where
 \begin{subequations}
\begin{align}
\begin{split}
\label{phi}
& \phi_\pm(\xi, \eta)  = {|\xi|}^2 - {|\xi-\eta|}^2 \pm |\eta| = 2\xi \cdot \eta -|\eta|^2 \pm |\eta| 
\end{split}
\\[.5 em]
\label{psi}
& \psi_{\pm}(\xi,\eta)  = \mp |\xi| - {|\xi-\eta|}^2 + {|\eta|}^2 = \mp |\xi| - {|\xi|}^2 + 2 \xi \cdot \eta \, .
\end{align}
\end{subequations}
From these  formulae, we can compute the space time resonance set for the system:
\begin{align*}
\mathscr{T_{\phi_\pm}} =& \{(\xi, \eta); \phi_\pm = 0\}, \,   
\\
\mathscr{S_{\phi_\pm}} =& \{(\xi, \eta); \nabla_\eta \phi_\pm = 0\} \, 
\\
\mathscr{R_{\phi_\pm}} =& \mathscr{T_{\phi_\pm}}\cap\mathscr{S_{\phi_\pm}} =   \{(\xi, \eta); \eta = 0, |\xi| =1/2 \} \, ,
\end{align*}
and
\begin{align*}
\mathscr{T_{\psi_\pm}} =& \{(\xi, \eta); \psi_\pm = 0\} \, , 
\\
\mathscr{S_{\psi_\pm}} =& \{(\xi, \eta); \nabla_\eta \psi_\pm = 0\} \, ,
\\
\mathscr{R_{\psi_\pm}} =& \mathscr{T_{\psi_\pm}}\cap\mathscr{S_{\psi_\pm}} =   \{(\xi, \eta); \xi = 0 \}  \, .
\end{align*}
With the space time resonance set identified, we proceed to look for null resonant interactions in the quadratic terms  \`a la \cite{PS}.  
For the profile $f$,  given by equation \eqref{inteqf},  the resonances are  null since  (see \cite{PS} page 4 and (3.2))
\begin{equation}
\label{dxiphi00}
\nabla_\xi \phi_\pm = - 2\eta = - 2 \frac{\eta}{|\eta|}(\frac{\eta}{|\eta|}\cdot \nabla_\eta \phi_\pm)-2\frac{\phi_\pm}{|\eta|}\frac{\eta}{|\eta|} \, .
\end{equation}
For the profile $g$ given by equation \eqref{inteqg} the resonances are  null since  
\begin{equation}
\label{symg0}
|\xi|=\frac{1}{2}\frac{\xi}{|\xi|} \cdot \nabla_\eta \psi_\pm \, .
\end{equation}
Since system \eqref{eq:profile} has null resonances, we can take advantage of them in combination with the space time resonance method.  
Roughly speaking we will obtain bounds on $u$ and $w$ in the following manner:
a) obtain standard energy bounds, namely $H^s$ bounds, on the profiles;
b) obtain weighted $L^2$ bounds on the profiles; and
c) obtain decay estimates  from Duhamel's formula or through the energy and weighted $L^2$ bounds.
The difficulties and ideas in proving these bounds will be illustrated below.

\subsubsection{Energy bounds} 
Energy estimates are  usually obtained in a straightforward fashion for semilinear equations. 
Because of the coupling $\Lambda|u|^2$ in equation \eqref{Z_0},  there is a derivative loss if one tries to obtain the estimates 
via Duhamel's formula or by multiplying by $(\partial_t u, \partial_t w_\pm)$.
This apparent derivative loss  can be handled by using a normal form transformation \cite{ShatahKGE} for high frequencies.

In addition to the derivative loss, we note here that good energy bounds can only be obtained if $w_{\pm}$ 
has an optimal decay of $t^{-1}$, due to the presence of $u w_\pm$ terms in the nonlinearity, 
whereas $u$ can be allowed to decay at a rate $t^{-(1+\alpha)}$ for some $0 < \alpha\le 1/2$.  
Bounds on high Sobolev norms are presented in section \ref{secenergy}.

\subsubsection{Weighted $L^2$ estimates for $G_\pm$}  
To obtain good weighted estimates one needs to use the  non resonant structure present in the system  \eqref{eq:profile}. 
Identity \eqref{symg0} allows us to estimate  $xG_\pm\in L^2$,  where 
 \begin{equation}
 \label{G}
G_\pm \overset{def}{=} \FF^{-1} \int_0^t \int_{\R^3}|\xi| e^{is \psi_{\pm}(\xi,\eta)}\hat f(\xi-\eta,s) \overline{\hat{f}} (\eta,s) d\eta ds \, .
 \end{equation}
Similarly we can obtain good bounds on $\Lambda {|x|}^2 G_\pm\in L^2$.   
These estimates are presented in section \ref{secweightedG}.   

The lack of good weighted bounds on ${|x|}^2 G_\pm$ leads to problems in obtaining weighted bounds on the nonlinearity in the \S equation
 \begin{equation}
 \label{F}
F_\pm \overset{def}{=}  \FF^{-1} \int_0^t \int_{\R^3} e^{is \phi_{\pm}(\xi,\eta)}\hat f(\xi-\eta,s) \hat g_\pm(\eta,s) d\eta ds \, ,
\end{equation}
as well as problems in obtaining the $L^\infty$ decay for $e^{it\Lambda} G_\pm$.  
Thus one has to find a new approach to obtain these bounds.

\subsubsection{$L^\infty$ bounds}  
The idea we utilize here consists in splitting $G_\pm$ according to the localization in real space of the input profiles\footnote{
A related idea is also exploited by Ionescu and Pausader in \cite{IP}, where the authors use 
norms based on a dyadic decomposition of the profiles in both frequency and real space.}  $f$.
Specifically, we split $f$ into a piece which is localized close to the origin and a far away piece,
by writing
\begin{equation*}
f_{\leq K} (x) = f(x) \rho \left( \frac{x}{K} \right) \quad , \quad  f_{\geq K} (x) = f(x) - f_{\leq K} (x) \, ,
\end{equation*}
for some smooth cutoff function $\rho$ with compact support which equals $1$ on the unit ball.
Then we split the profile $f$ in \eqref{G} into a localized piece  $f_{\leq s^{1/8}}$, and a far away piece $f_{\geq s^{1/8}}$.
This allows us to gain decay on the pieces of $G_\pm$
\begin{align*}
\hat{G}_1 & \overset{def}{=}   \int_0^t \int_{\R^3} |\xi| e^{is \psi_{\pm}(\xi,\eta)} \hat{f_{\leq s^{1/8} }} (\xi-\eta,s) \overline{ \hat{f} }(\eta,s) d\eta ds
\\
\hat{G}_2 &  \overset{def}{=}   \int_0^t \int_{\R^3} |\xi| e^{is \psi_{\pm}(\xi,\eta)} \hat{f_{\geq s^{1/8} }} (\xi-\eta,s) \overline{ \hat{f_{\geq s^{1/8} }} } (\eta,s) d\eta ds \, ,
\end{align*}
by  noting the following:

\begin{enumerate}

\item For the localized component, i.e., $f_{ \leq s^{1/8}}$, one has that
  $\|e^{it\Delta} f_{\leq s^{1/8}}\|_{L^\infty}$ decays faster than ${\| e^{is\Delta} f \|_{L^\infty}}$ does,
  and ${\| {|x|}^2 f_{\leq s^{1/8}} \|_{L^2}}$ grows slower than ${\| {|x|}^2 f \|_{L^2}}$ does, see \eqref{apriorif}.
  These facts, plus the non resonance structure, allow us to show essentially $\Lambda^2 G_1 \in L^1$, which is sufficient to obtain  
  the $L^\infty$ decay of $e^{\mp it\Lambda} G_1$.

\item For the far away component $f_{\geq s^{1/8}}$, 
	one has $\| f_{ \geq s^{1/8}}\|_{L^2} \lesssim s^{-1/8} \|xf\|_{L^2}$. 
	This fact, and again the non resonance structure, allow us to gain time decay for the integrand in $G_2$.
	Combining these and the dispersive estimate for the wave operator we obtain the $L^\infty$ decay of  $e^{\mp it\Lambda} G_2$.
\end{enumerate}
The decay estimate on $e^{\mp it\Lambda} G_\pm$ is presented in section \ref{secdecayw}.

\subsubsection{Refined estimates on $G_\pm$}  
Since  $\nabla_\xi\psi_\pm$ does not vanish on $\mathscr{R}_{\psi_\pm}$,  we can only obtain lousy bounds on $|x|^2G_\pm$ in $L^2$. 
To deal with this difficulty, we again split $G_\pm$ into two parts: 
a component $g_1$ which comes from localized interactions, and $g_2$ which has at least one term far away
\begin{align*}
\hat{g}_1 & \overset{def}{=}  \int_0^t \int_{\R^3} |\xi| e^{is \psi_{\pm}(\xi,\eta)} \hat{f_{\leq s^{1/4} }} (\xi-\eta,s) \overline{ \hat{f_{\leq s^{1/4} }} }(\eta,s) d\eta ds
\\
\hat{g}_2 & \overset{def}{=} \int_0^t \int_{\R^3} |\xi| e^{is \psi_{\pm}(\xi,\eta)} \hat{f_{\geq s^{1/4} }} (\xi-\eta,s) \overline{ \hat{f} } (\eta,s) d\eta ds \, .
\end{align*}
The term $g_1$ has well localized inputs, and thus we have good estimates on $|x|^2 g_1\in L^2$, see Lemma \ref{Lemmag_1}. 
The term $g_2$ instead has good small frequency behavior as shown in Lemma \ref{Lemmag_2}.

\subsubsection{Weighted $L^2$ estimates for $F_\pm$}  
Estimates of $|x|F\in L^2$ can be obtained in a relatively straightforward manner using \eqref{dxiphi00}.  
To estimate $|x|^2F\in L^2$ we use the splitting of $G_\pm$ above, and see that we need to control bilinears term of the form 
\begin{equation*}
B(f,g_i) (t,\xi) = \int_0^t \int_{\R^3} s^2 \eta^2  e^{is \phi_{\pm}(\xi,\eta)} \hat f(\xi-\eta) \hat g_i (\eta) d\eta ds \, , \quad  i=1,2 \, .
\end{equation*}
The term $B(f,g_1)$ can be estimated by integrating by parts twice in $\eta$, again via \eqref{dxiphi00},
since $|x|^2 g_1$ does not grow too fast.  
For the term involving  $g_2$,  we use the fact that $g_2$ has good behavior for small frequencies 
to  excise a relatively large neighborhood of $\eta = 0$ around the space time resonant set.
We can then control $B(f,g_2)$ in this neighborhood. 
On the complement we can integrate by parts in $\eta$ twice and use the available bound on the $L^2$ norm of $\Lambda |x|^2 g_2$ 
to eventually control $B(f,g_2)$.   
These estimates are presented in section \ref{secweightedF}.

\section{Norms and bounds}\label{secnorms}
Our proof of Theorem  \ref{maintheo} consists of closing a bootstrap argument with the following a priori bounds:
\begin{equation}
\label{apriorif}
 \quad  {\| f(t) \|}_{H^{N+1}}  \lesssim \e_0 t^\d  \, , 
\quad
{\| x f(t) \|}_{L^2} \lesssim \e_0 t^\d \, , 
\quad{\| {| x |}^2 f(t) \|}_{L^2} \lesssim \e_0 t^{1 - 2\a - \d}  \, , 
\end{equation}
and
\begin{align}
\label{apriorig}
& {\|  g_\pm(t) \|}_{H^N}  \lesssim \e_0 \, ,
\quad
{\| e^{\mp it\Lambda} g_\pm(t) \|}_{ \dot{B}^0_{\infty,1} } \lesssim \frac{\e_0}{t} \, , 
\\
\label{aprioriG}
& {\| x G_\pm(t) \|}_{H^1} \lesssim \e_0  \, , 
\quad
{\| \Lambda {| x |}^2  G_\pm(t) \|}_{L^2} \lesssim \e_0 t^{ \beta }   \, ,
\end{align} 
where the parameters are chosen such that
\begin{align}
\label{param}
\a = \frac{1}{6} - 2 \d
\quad , \qquad  
\b = 1 - 3 \a
\quad , \qquad
\frac{5}{N} \leq \d \quad , \qquad  \d \ll 1 \, .
\end{align}
Here and in what follows, we denote by $\dot{B}^s_{p, q}$ the Besov space defined by the norm 
\begin{equation}\label{Besov}
\|u\|_{\dot{B}^s_{p,q}}:=\left\| 2^{sk}\|P_k u\|_{L_x^p(\R^3)}\right\|_{l_k^q(\Z)}
\end{equation}
where $P_k$ denotes the Littlewood-Paley projection onto frequencies $|\xi|\sim 2^k$.

We define the norm $X$ associated to the bounds \eqref{apriorif} and \eqref{apriorig} by
\begin{align}
\nonumber
{\|(u,w_\pm)\|}_{X}  \stackrel{def}{=} \sup_t   \left( t^{-\d}  {\|  f(t)  \|}_{H^{N+1}} \right. & +  t^{-\d}{\| xf (t) \|}_{L^2}  
    +  t^{-1+2\a+\d}  {\| {|x|}^2 f (t) \|}_{L^2}
\\
\label{norm}
 & + \, {\|  g_\pm(t) \|}_{H^{N}} + t \left. {\|  e^{\mp it\Lambda} g_\pm (t) \|}_{\dot{B}^0_{\infty,1} }  \right) \, .
\end{align}
From \eqref{inteqf}--\eqref{inteqg}, and the definition of $F_\pm$ and $G_\pm$ in \eqref{F} and \eqref{G},
we have 
\begin{align*}
f (t,x) & = f_0 (x) \mp \sum_{\pm} i F_\pm (t,x)
\\
g_\pm (t,x) & = g_\pm (0,x) - i G_\pm (t,x) \, ,
\end{align*}
where $f_0(x) = u_0(x)$ and $ g_\pm (0,x) = \pm n_0(x) + i \Lambda^{-1} n_1(x)$.
From the hypotheses \eqref{data} on the initial data we have (see also \eqref{decaywave0} below)
\begin{align*}
{\left\| \left( e^{it\Delta} f_0, e^{\mp it\Lambda} g_\pm (0) \right) \right\| }_X \leq \e_0 \, .
\end{align*}
In Propositions \ref{proenergy}, \ref{prodecayw} and \ref{proweightedF}, we will show that
\begin{align*}
{\left\| \left( e^{it\Delta} F_\pm (x), e^{\mp it\Lambda} G_\pm \right) \right\| }_X \lesssim {\|(u,w_\pm)\|}^2_{X} \,  
\end{align*}
provided ${\|(u,w_\pm)\|}_{X}$ is small enough, and this will imply
\begin{align*}
{\left\| ( u, w_\pm ) \right\| }_X \leq \e_0 + C {\|(u,w_\pm)\|}^2_{X} \, .
\end{align*}
A standard continuation argument will then guarantee a global solution in the space defined by the norm \eqref{norm},
provided this is small enough.

We remark here that the weighted bounds \eqref{aprioriG} on $G_\pm$, which are obtained in Proposition \ref{proweightedG},
are only instrumental to the proof of the weigthed bounds on $F_\pm$ given in Proposition \ref{proweightedF}.
We do not have estimates like \eqref{aprioriG} for $g_\pm$,
as these would require some vanishing moment condition on the data $n_1$.

\begin{remark}[Linear dispersive estimates]
Note that from the linear estimates for the Schr\"odinger group
\begin{align}\label{disp}
& {\|e^{it\Delta} f\|}_{L^6} \lesssim  \frac{1}{t} {\| x f \|}_{L^2}
\quad , \quad 
  {\| e^{it\Delta} f \|}_{L^\infty} \lesssim  \frac{1}{t^{\frac 32}} {\| x f \|}^{\frac 12}_{L^2} {\| x^2 f \|}^{\frac 12}_{L^2} \, ,
\end{align}
we deduce  that the $X$ norm bounds
\begin{align}
\label{SL^6}
{\| e^{it\Delta} f \|}_{L^6}  & \lesssim \frac{1}{ { t }^{1-\d}} {\| u \|}_X \, ,
\\
\label{decayu}
{\| e^{it\Delta} f \|}_{L^\infty} & \lesssim \frac{1}{ { t }^{1+\a}} {\| u \|}_X  \, .
\end{align}
Moreover, by the linear dispersive estimate for the wave equation
\begin{equation}
\label{linearwave0}
{\| e^{is\Lambda} h \|}_{\dot{B}^0_{p,r}} 
				\lesssim \frac{1}{t^{1-\frac{2}{p}}} {\| h \|}_{ {\dot{B}^{2(1-2/p)}_{{p^\prime},r}} } \quad , \quad p \geq 2 \, ,
\end{equation}
(cf. for example \cite{SS}), and the fact that $g_\pm(0) = \Lambda^{-1} i n_1 \pm n_0$, 
we see that \eqref{data2} implies
\begin{equation}
\label{decaywave0}
{\| e^{\mp i t\Lambda} g_{\pm}(0) \|}_{ \dot{B}^0_{\infty,1} } \lesssim \frac{\e_0}{t} \, .
\end{equation}
Finally, we note that by \eqref{linearwave0} with $r = 2$, and embeddings between Besov and Sobolev spaces, we have
\begin{equation}
\label{linearwave}
{\| e^{is\Lambda} h \|}_{L^p} 
				\lesssim \frac{1}{ t^{1-2/p} } {\left\|  \Lambda^{2(1-2/p)} h \right\|}_{ L^{p^\prime} } \, .
\end{equation}

\end{remark}

\vskip5pt
\begin{convention}
The cases $\pm$ will be treated identically in our analysis. Therefore for ease of exposition we will  drop the apex $\pm$.
\end{convention}

\vskip10pt
\section{Energy Estimates}\label{secenergy}
In this section we are going to prove the following:
\begin{proposition}
\label{proenergy}
Let $G$ and $F$ be given by \eqref{G} and \eqref{F} respectively. Then, for ${\| (u, w) \|}_X \lesssim \e_0$ we have
\begin{equation*}
 {\| G \|}_{H^{N}}  + t^{-\d} {\| F \|}_{H^{N+1}} \lesssim {\| (u,w) \|}_X^2 \, .
\end{equation*}
\end{proposition}

\begin{convention}
As the norm ${\| (u,w) \|}_X^2$ will appear at the end of all our chains of inequalities, we omit it for lighter notations.
\end{convention}

\noindent
To prove Proposition \eqref{proenergy} we start with the easy estimate for the $H^N$ norm of the acoustic wave component $G$.
Then, in order to control the $H^{N+1}$ norm of $F$ we will use a normal form transformation to make up for a derivative loss.
We will then use the control on these Sobolev norms to reduce all of our estimates to frequencies smaller than $s^{\delta_N}$, 
where $\delta_N \ll 1$ is chosen small depending on $N$.

\subsection{Estimate  ${\| G \|}_{H^N} \lesssim 1$.}
This bound follows just by H\"{o}lder's inequality and the more than integrable decay of the \S component $u$:
\begin{align*}
{\| G \|}_{H^N} \lesssim \int_0^t {\left\| \Lambda {| u |}^2 \right\|}_{H^N} \, ds  \lesssim \int_0^t {\| u \|}_{H^{N+1}} {\| u \|}_{L^\infty} \, ds
						\lesssim \int_0^t s^\d \frac{1}{ {\l s \r}^{1+\a}} \, ds \lesssim 1 \, .
\end{align*}
Here we used \eqref{decayu} for times $s \geq 1$, and Sobolev embedding for $s \leq 1$ to deduce that
${\| u \|}_{L^\infty} \lesssim  {\l s \r}^{ - 1 - \a}$.

\subsection{Estimate ${\| F \|}_{H^{N+1}} \lesssim t^\d$}  First note that by \eqref{inteqf} and since ${\| f \|}_{L^2}$ is conserved,  
we only need to estimate the $\dot{H}^{N+1}$ norm of $F$ in \eqref{F}.
Let us define smooth positive radial cutoff functions $\chi_1(\xi,\eta)$ and $\chi_2(\xi,\eta)$, with $\chi_1 + \chi_2 = 1$,
and such that 
\[
\chi_2 (\xi,\eta) = 1 \, \text{ if } \, 100 |\xi -\eta| \leq |\eta| \, , 
\, \text{ and } \, \chi_2 (\xi,\eta) = 0  \, \text{ if } \, |\eta| \leq 50 |\xi-\eta| \, .
\]
We then write $F = F_1 + F_2$ where
\begin{align}
\label{F_1}
F_1 (t,x) & := \FF^{-1} \int_0^t \int_{\R^3} \chi_1(\xi,\eta) e^{is \phi (\xi,\eta)}\hat{f}(\xi-\eta,s) \hat{g}(\eta,s) d\eta ds
\\
\label{F_2}
F_2 (t,x) & := \FF^{-1} \int_0^t \int_{\R^3} \chi_2(\xi,\eta) e^{is \phi (\xi,\eta)}\hat{f}(\xi-\eta,s) \hat{g}(\eta,s) d\eta ds \, .
\end{align}
On the support of $\chi_1$ we have $|\xi -\eta| \gtrsim |\eta|$ and hence derivatives applied to $F_1$ fall only on $u = e^{is\Delta} f$.
This term is then easily estimated using the Coifmain-Meyer theorem \cite{CM}:
\begin{align*}
{\| F_1 \|}_{H^{N+1}} \lesssim \int_0^t  {\| u \|}_{H^{N+1}} {\| w \|}_{L^\infty}  \, ds
						\lesssim \int_0^t s^\d \frac{1}{ s } \, ds \lesssim t^\d \, .
\end{align*}
In estimating $F_2$ we can reduce ourselves to the case $|\eta| \geq 100$, otherwise an application of H\"{o}lder's inequality as above suffices.
Then we observe that on the suppport of $\chi_2$, the phase $\phi_\pm$ in \eqref{phi} satisfies the following:
\begin{equation}
\left| \phi_\pm(\xi,\eta) \right| \geq {|\eta|}^2 - |\eta| - 2|\xi-\eta| |\eta| \gtrsim {|\eta|}^2 \, .
\end{equation}
This lower bound will allow us to recover the loss of derivative that would occur when all $N+1$ derivatives fall on $w = e^{it\Lambda} g$.
We have 
\begin{align*}
\Lambda^{N+1} F_2 (t,x) & = \FF^{-1} \int_0^t \int_{\R^3} {|\xi|}^{N+1} \chi_2(\xi,\eta)  
      e^{is \phi (\xi,\eta)}\hat{f}(\xi-\eta,s) \hat{g}(\eta,s) d\eta ds \, .
\end{align*}
We then integrate by parts in $s$ to get
\begin{subequations}
\begin{align}
\label{F_21}
\Lambda^{N+1} F_2 (t,x) & = {\left. \FF^{-1} \int_{\R^3} \frac{\chi_2(\xi,\eta) {|\xi|}^{N+1}}{ i {|\eta|}^N \phi (\xi,\eta)}  
						 e^{is \phi (\xi,\eta)} \hat{f}(\xi-\eta,s)  {|\eta|}^{N} \hat{g} (\eta,s) d\eta
	    \right|}_{s=0}^{s=t}
\\
\label{F_22}
& - \FF^{-1} \int_0^t \int_{\R^3}\frac{\chi_2(\xi,\eta) {|\xi|}^{N+1}}{ i {|\eta|}^N \phi(\xi,\eta)}  
					e^{is \phi (\xi,\eta)}  \partial_s \hat{f} (\xi-\eta,s) {|\eta|}^{N} \hat{g} (\eta,s) d\eta ds
\\
\label{F_23}
& - \FF^{-1} \int_0^t \int_{\R^3} \frac{\chi_2(\xi,\eta) {|\xi|}^{N+1}}{ i {|\eta|}^N \phi (\xi,\eta)}
					e^{is \phi (\xi,\eta)} \hat{f} (\xi-\eta,s)   {|\eta|}^{N} \partial_s \hat{g} (\eta,s) d\eta ds \, .
\end{align}
\end{subequations}
Thanks to the lower bound on $\phi_\pm$ one easily verifies that  on the support of $\chi_2$
\begin{equation}
\left| \frac{\chi_2(\xi,\eta)  {|\xi|}^{N+1} }{ {|\eta|}^N \phi_\pm(\xi,\eta)}  \right| \lesssim 1 \, .
\end{equation}
Disregarding the contribution from $s=0$ which is easier to estimate, 
the term \eqref{F_21} above is estimated by Plancharel's theorem and Young's inequality:
\begin{align*}
{\| \eqref{F_21} \|}_{L^2} & \lesssim {\left\| \int_{\R^3}  \frac{\chi_2(\xi,\eta)  {|\xi|}^{N+1} }{ {|\eta|}^N \phi_\pm(\xi,\eta)}
						 e^{it \phi_{\pm}(\xi,\eta)} \hat{f}(\xi-\eta,t)  {|\eta|}^{N} \hat{g}_\pm (\eta,t) d\eta \right\|}_{L_\xi^2}
\\
& \lesssim {\left\| \int_{\R^3}  \left| \hat{u}(\xi-\eta,t) \right|  \left| {|\eta|}^{N} \hat{w}_\pm (\eta,t) \right| d\eta \right\|}_{L_\xi^2}
\lesssim {\| u \|}_{H^2} {\| w \|}_{H^N}  \lesssim t^\d \, .
\end{align*}
Using a similar argument, and the fact that $e^{is\Delta} \partial_s f = u w$, we can bound
\begin{align*}
{\| \eqref{F_22} \|}_{L^2} & \lesssim \int_0^t {\left\| e^{is\Delta} \partial_s f  \right\|}_{H^2} {\| w \|}_{H^N}  \, ds 
\\
& \lesssim \int_0^t \left[ {\| w \|}_{H^2} {\| u \|}_{L^\infty} + {\| u \|}_{H^2} {\| w \|}_{L^\infty} \right]  {\| w \|}_{H^N}  \, ds 
		\lesssim \int_0^t \frac{1}{ {\l s \r}^{1+\a}} + s^\d \frac{1}{s} \, ds \lesssim t^\delta \, .
\end{align*}
Analogously, since $e^{is\Lambda} \partial_s w = \Lambda {|u|}^2 $, we have
\begin{align*}
{\| \eqref{F_23} \|}_{L^2} & \lesssim \int_0^t {\| u  \|}_{H^2} {\left\| e^{is\Lambda} \partial_s w  \right\|}_{H^N}  \, ds 
		\\
		& \lesssim \int_0^t {\| u \|}_{H^2} {\| u \|}_{L^\infty}  {\| u \|}_{H^{N+1}}  \, ds 
		\lesssim \int_0^t s^\d \frac{1}{ {\l s \r}^{1+\a}} s^\d \, ds \lesssim t^\delta \, 
\end{align*}
since $\alpha > \delta$.

\subsection{High frequency cutoff}
In the previous section we have established the a priori bounds ${\| u \|}_{H^{N+1}} \lesssim t^\d$ and ${\| n \|}_{H^N} \lesssim 1$.
Let us denote by $P_{\geq k}$ the Littlewood-Paley projection on frequencies larger or equal to $2^k$.
Since for $k \geq 0$ one has
\begin{align*}
{\| P_{\geq k} v (s) \|}_{L^2} & \lesssim 2^{-k l } {\| v(s) \|}_{H^l} \, ,
\end{align*}
then, for frequencies $2^k \gtrsim s^{ 2/(N-2) }$ we have
\begin{align}
\label{Pk1}
{\| P_{\geq k} u (s) \|}_{H^3} & \lesssim  2^{-k (N-2) } {\| u(s) \|}_{H^{N+1}} \lesssim \frac{1}{ {\l s \r}^2} s^\d
\\
\label{Pk2}
{\| P_{\geq k} w (s) \|}_{H^2} & \lesssim  2^{-k (N-2) } {\| w(s) \|}_{H^{N}} \lesssim \frac{1}{ {\l s \r}^2} \, .
\end{align}
This shows that in estimating weighted norms of the bilinear terms $F$ and $G$ in \eqref{F} and \eqref{G},
we can  always reduce our analysis to frequencies $|\xi-\eta| , |\eta| \lesssim s^{2/(N-2)}$.
Indeed, if at least one of the frequencies $|\eta|$ or $|\xi-\eta|$ is greater than $s^{2/(N-2)}$,
all the desired bounds can be shown to hold true in a straightforward fashion.
This is because of \eqref{Pk1} and \eqref{Pk2} above, 
and because the action of weights on  Littlewood-Paley projections $P_{\geq k}$, with $2^k \gtrsim s^{2/(N-2)}$,
is harmless, and would only give contributions which are much easier to treat than the ones we are going to estimate below.
Therefore we agree on the following:

\noindent
\begin{convention}
In the rest of the paper, we assume that all frequencies $|\xi-\eta|$ and $|\eta|$ appearing in the 
estimates of the bilinear terms \eqref{F} and \eqref{G}, are bounded above by $s^{\delta_N}$ where  $\delta_N := \frac{2}{N-2}$, 
and the integer $N \gg 1$ is determined in the course of our proof by several upperbounds on $\d_N$.
In particular, expressions such as $|\xi|$ or $\nabla_\xi \psi (\xi,\eta)$ will be often replaced by a factor of $s^{\delta_N}$.
\end{convention}

\section{Weighted estimates for the wave component}\label{secweightedG}
In this section we show the bounds on weighted $L^2$ norms of $G$. 
\begin{proposition}
\label{proweightedG}
Let $G$ be the bilinear term defined in \eqref{G}, then 
\begin{equation*}
{\| x G \|}_{H^1} + t^{-1+3\a} {\| \Lambda x^2 G \|}_{L^2}  \lesssim {\| (u,w) \|}_X^2 \, .
\end{equation*}
\end{proposition}
\noindent
To prove this it is crucial to notice the presence of a null resonant structure in the nonlinear term $G$, 
that is to say, the vanishing of the symbol on the space resonant set.
More precisely we have
\begin{equation}
\label{symg}
|\xi|=\frac{1}{2}\frac{\xi}{|\xi|} \cdot \nabla_\eta \psi \, ,
\end{equation}
which allows us to integrate by parts gaining decay in $s$.

\subsection{Estimate ${\| x G \|}_{H^1} \lesssim 1$}
Applying $\nabla_\xi$ to $\hat{g}$ gives the terms:
\begin{align}
\label{dxig1}
& \int_0^t \int_{\R^3} e^{i \psi (\xi,\eta)s} |\xi| \nabla_\xi \hat{f}(\xi-\eta,s) \overline{\hat{f}} (\eta,s) d\eta ds
\\
\label{dxig2}
& \int_0^t \int_{\R^3} s \nabla_\xi \psi  e^{i s \psi (\xi,\eta)} |\xi| \hat{f}(\xi-\eta,s) \overline{\hat{f}} (\eta,s) d\eta ds \, ,
\end{align}
plus an easier term when $\nabla_\xi$ hits the symbol $|\xi|$.
\eqref{dxig1} is easily estimated by H\"{o}lder's inequality:
\begin{align*}
{\| \eqref{dxig1} \|}_{L^2} & \lesssim \int_1^t s^{\delta_N} {\|  x f \|}_{L^2} {\|  e^{is\Delta} f  \|}_{L^\infty} \, ds
\lesssim \int_0^t s^{\delta_N} s^\d \frac{1}{ {\langle s \rangle}^{1+\a}} \, ds \lesssim 1 
\end{align*}
since $\alpha>\delta+\delta_N$.
Using the identity \eqref{symg}
and integrating by parts in $\eta$, \eqref{dxig2} gives terms like \eqref{dxig1}. Therefore we can skip them.

\vskip10pt
\subsection{Estimate ${\| \Lambda x^2 G \|}_{L^2} \lesssim t^{1 - 3\alpha}$}\label{lamdax2G}
Applying $|\xi| \nabla_\xi^2$ to $\hat{g}$ gives the following main\footnote{
The remaining terms where at least one derivative $\nabla_\xi$ hits the symbol $|\xi|$ are easier to estimate.} contributions:
\begin{subequations}
\begin{align}
\label{xidxi^2g1}
& \int_0^t \int_{\R^3} e^{i s \psi (\xi,\eta)} {|\xi|}^2 \nabla_\xi^2 \hat{f}(\xi-\eta,s) \overline{\hat{f}} (\eta,s) d\eta ds
\\
\label{xidxi^2g2}
& \int_0^t \int_{\R^3} s \nabla_\xi \psi  e^{i s \psi (\xi,\eta)} {|\xi|}^2  \nabla_\xi \hat{f}(\xi-\eta,s) \overline{\hat{f}} (\eta,s) d\eta ds
\\
\label{xidxi^2g3}
& \int_0^t \int_{\R^3} s^2 {(\nabla_\xi \psi)}^2 e^{i s \psi (\xi,\eta)} {|\xi|}^2 \hat{f}(\xi-\eta,s) \overline{\hat{f}} (\eta,s) d\eta ds \, .
\end{align}
\end{subequations}
\eqref{xidxi^2g1} can be directly estimated as follows:
\begin{align*}
{\| \eqref{xidxi^2g1} \|}_{L^2} & \lesssim \int_0^t s^{2\delta_N} {\|  x^2 f \|}_{L^2} {\|  e^{is\Delta} f  \|}_{L^\infty} \, ds
\\
& \lesssim \int_0^t s^{2\delta_N} s^{1-2\a-\d} \frac{1}{s^{1+\a}} \, ds \lesssim t^{1 - 3\a - \d + 2 \delta_N} \, ;
\end{align*}
this satisfies the desired bound provided $2\delta_N \leq \d$.

Using the identity \eqref{symg} and integrating by parts in $\eta$ once in \eqref{xidxi^2g2}, and twice in \eqref{xidxi^2g3},
gives terms similar to \eqref{xidxi^2g1}, plus the following:
\begin{align}
\label{xidxi^2g4}
& \int_0^t \int_{\R^3} e^{i s \psi (\xi,\eta)} m_2 (\xi,\eta)  \nabla_\eta \hat{f}(\xi-\eta,s) \nabla_\eta \overline{\hat{f}} (\eta,s) d\eta ds \, ,
\end{align}
where $m_2$ is a symbol satisfying homogeneous bounds of order $2$ for large frequencies,
and is otherwise harmless.
Using the dispersive estimate and ${\| \cdot \|}_{L^{4/3}} \lesssim {\|  \cdot \|}_{L^2}^{1/4} {\| x \cdot \|}_{L^2}^{3/4}$,
this can be bounded as follows
\begin{align*}
{\| \eqref{xidxi^2g4} \|}_{L^2} & \lesssim  \int_0^t s^{2\delta_N} {\|  e^{is\Delta} x f \|}^2_{L^4} \, ds
		\lesssim \int_0^t s^{2\delta_N} \frac{1}{ s^{3/2} }{\|  x f \|}_{L^{4/3}}^2 \, ds
\\
& \lesssim \int_0^t s^{2\delta_N} \frac{1}{ s^{3/2} }  {\|  x^2 f \|}_{L^2}^{3/2} {\| x f \|}_{L^2}^{1/2}  \, ds
		\lesssim \int_0^t s^{2\delta_N} \frac{1}{ s^{3/2} } s^{ \frac{3}{2} (1 -2\alpha -\d) } s^{\d/2} \, ds 
		\lesssim t^{1-3\alpha}
\end{align*}
again provided $2 \delta_N \leq \d$.

\section{Decay estimate for the wave component}
\label{secdecayw}
From \eqref{decaywave0} we already have the necessary pointwise decay for $e^{it\Lambda} g_\pm(0)$.
We then need to show the following:
\begin{proposition}
\label{prodecayw}
Let $G$ be the bilinear term defined in \eqref{G}, then 
\begin{equation*}
t {\| e^{it\Lambda} G \|}_{ \dot{B}^0_{\infty,1} } \lesssim {\| (u,w) \|}_X^2 \, .
\end{equation*}
\end{proposition}

\noindent
In order to prove the above Proposition we split $G$ into two parts, depending on the localization of the inputs.
More precisely, we let
$G = G_1 + G_2$
where 
\begin{align*}
G_1 & :=   G(f_{\leq s^{1/8} } , \bar{f}) + G(f_{\geq s^{1/8}}, \bar{f}_{\leq s^{1/8} })
\\
G_2 & :=   G(f_{\geq s^{1/8} } , \bar{f}_{\geq s^{1/8} } ) \, .
\end{align*}
The component $G_1$ will be shown to be bounded in $\dot{B}^{2}_{1,1}$, which gives the desired bound on $e^{it\Lambda} G_1$.
The decay of $e^{it\Lambda} G_2$ will instead be proven using the null structure \eqref{symg}, and the fact that the
$L^2$ norm of $f_{\geq s^{1/8}}$ decays in $L^2$.

Since the statement of Proposition \ref{prodecayw} is easy to obtain for $t \leq 1$, 
in the integral \eqref{G} which defines $G$ we will only consider the contribution going from $1$ to $t$.
Also, since the two terms in the definition of $G_1$ are similar, we can reduce to consider $G_1$ and $G_2$ given by
\begin{align}
\label{G_1}
\hat{G}_1 & =   \int_1^t \int_{\R^3} |\xi| e^{is \psi (\xi,\eta)} \hat{f_{\leq s^{1/8} }} (\xi-\eta,s) \overline{ \hat{f} }(\eta,s) d\eta ds
\\
\label{G_2}
\hat{G}_2 & =   \int_1^t \int_{\R^3} |\xi| e^{is \psi (\xi,\eta)} \hat{f_{\geq s^{1/8} }} (\xi-\eta,s) \overline{ \hat{f_{\geq s^{1/8} }} } (\eta,s) d\eta ds \, .
\end{align}

\subsection{Decay estimate for $e^{it\Lambda} G_1$}
To show that $G_1$ is bounded in  $\dot{B}^{2}_{1,1}$ we will interpolate weighted $L^2$ norms inside the time integral.
This type of argument was also used in \cite{PS}.
Here we will crucially use the ``small'' support of $f_{\leq s^{1/8}}$ to get improvements on its weighted norms, 
and on the decay of $e^{is\Delta}f_{\leq s^{1/8}}$.
Recalling that we are only considering frequencies $k$ such that $2^k \leq s^{\d_N}$,
we aim to prove
\begin{align*}
\int_1^t  \, \sum_{k = -\infty}^{ \log s^{\d_N} } 2^{2k} {\| P_k \Lambda e^{-is\Lambda} 
				\left( e^{is\Delta} f_{\leq s^{1/8}} e^{-is\Delta} \bar{f} \right) \|}_{L^1} \, ds  \lesssim 1 \, .
\end{align*}
Converting a factor of $2^k$ into a derivative $\Lambda$, throwing away the projection $P_k$, and performing the sum,
we see that is suffices to show
\begin{align*}
\int_1^t  s^{\d_N} {\| \Lambda^2 e^{-is\Lambda} \left( e^{is\Delta} f_{\leq s^{1/8}} e^{-is\Delta} \bar{f} \right) \|}_{L^1} \, ds  \lesssim 1 \, .
\end{align*}
Since ${\| \cdot \|}_{L^1} \lesssim {\| x \cdot \|}_{L^2}^{1/2} {\| x^2 \cdot \|}_{L^2}^{1/2}$, 
the above estimate will follow from the inequalities
\begin{subequations}
\begin{align}
\label{decayg_11}
& {\left\| |x| \Lambda^2 e^{-is\Lambda} \left( e^{is\Delta} f_{\leq s^{1/8}} e^{-is\Delta} \bar{f} \right) \right\|}_{L^2}  \lesssim  \frac{1}{s^{7/4}} \, ,
\\
\label{decayg_12}
& {\left\| {|x|}^2 \Lambda^2 e^{-is\Lambda} \left( e^{is\Delta} f_{\leq s^{1/8}} e^{-is\Delta} \bar{f} \right) \right\|}_{L^2} \lesssim \frac{1}{\sqrt{s}} \, .
\end{align}
\end{subequations}
The quantity measured in the $L^2$ norm in \eqref{decayg_11} is given by a sum of terms of the form
\begin{subequations}
\begin{align}
\label{xg_1}
& \int_{\R^3} {|\xi|}^2 e^{is \psi (\xi,\eta)} \nabla_\xi \hat{f_{\leq s^{1/8} }} (\xi-\eta,s) \overline{ \hat{f} }(\eta,s) d\eta \, ,
\\
\label{xg_2}
& \int_{\R^3} {|\xi|}^2 e^{is \psi (\xi,\eta)} \,s \nabla_\xi \psi   \hat{f_{\leq s^{1/8} }} (\xi-\eta,s) \overline{ \hat{f} }(\eta,s) d\eta \, ,
\end{align}
\end{subequations}
plus similar or easier ones.
Using one of the factors $|\xi|$ in \eqref{xg_1} and two of them in \eqref{xg_2}, we integrate by parts obtaining 
as main contributions
\begin{subequations}
\begin{align}
\label{xg_11}
& \frac{1}{s} \int_{\R^3} m_1(\xi,\eta) e^{is \psi  (\xi,\eta)} \nabla^2_\eta \hat{f_{\leq s^{1/8} }} (\xi-\eta,s) \overline{ \hat{f} }(\eta,s) d\eta  \, ,
\\
\label{xg_12}
& \frac{1}{s} \int_{\R^3} m_1(\xi,\eta) e^{is \psi  (\xi,\eta)} 
			\nabla_\eta \hat{f_{\leq s^{1/8} }} (\xi-\eta,s) \nabla_\eta \overline{ \hat{f} }(\eta,s) d\eta  \, ,
\\
\label{xg_21}
& \frac{1}{s} \int_{\R^3} m_1(\xi,\eta) e^{is \psi  (\xi,\eta)} \hat{f_{\leq s^{1/8} }} (\xi-\eta,s) \nabla_\eta^2 \overline{ \hat{f} }(\eta,s) d\eta \, ,
\end{align}
\end{subequations}
where $m_1(\xi,\eta)$ denotes a symbols with homogenous bounds of order $1$ for large frequencies
and which is otherwise  harmless. 
%
%
%
The first term is then estimated as follows:
\begin{align*}
{\| \eqref{xg_11} \|}_{L^2} & \lesssim 
\frac{1}{s} s^{\delta_N} {\|  x^2 f_{\leq s^{1/8}} \|}_{L^2} {\| e^{is\Delta}  f  \|}_{L^\infty}
\\
& \lesssim \frac{1}{s} s^{\delta_N} s^{1/8} s^\d  \frac{1}{s^{1+\a}}  \lesssim  \frac{1}{s^{7/4}}
\end{align*}
where we used \eqref{param} with $\d$ small enough. 
Similarly we can bound
\begin{align*}
{\| \eqref{xg_12} \|}_{L^2} & \lesssim 
\frac{1}{s} s^{\delta_N} {\|  e^{is\Delta} x f_{\leq s^{1/8}} \|}_{L^\infty} {\|  x f  \|}_{L^2}
\\
& \lesssim \frac{1}{s} s^{\delta_N} \frac{1}{ s^{3/2} } {\| x f_{\leq s^{1/8}} \|}_{L^1}  {\|  x f  \|}_{L^2}
\\
& \lesssim \frac{1}{s} s^{\delta_N} \frac{1}{ s^{3/2} } s^\d  s^{3/16} s^\d
\lesssim  \frac{1}{s^{7/4}}
\end{align*}
Finally we can estimate
\begin{align*}
{\| \eqref{xg_21} \|}_{L^2} & \lesssim 
\frac{1}{s} s^{\delta_N} {\|  e^{is\Delta} f_{\leq s^{1/8}} \|}_{L^\infty} {\| x^2 f  \|}_{L^2}
\\
& \lesssim \frac{1}{s} s^{\delta_N} \frac{1}{s^{3/2}}  {\| f_{\leq s^{1/8}} \|}_{L^1}  s^{1-2\a-\d}
\\
& \lesssim \frac{1}{s} s^{\delta_N} \frac{1}{s^{3/2}} s^{1/16}  s^\d s^{1-2\a-\d}
\lesssim \frac{1}{s^{7/4}}
\end{align*}
where in the last inequality we have used $1- 2\a +\d_N = \frac{2}{3} + 4\d + \d_N \leq \frac{11}{16}$,
in accordance with \eqref{param} for $\d$ small enough.

Having established \eqref{decayg_11}, it is then easy to see that \eqref{decayg_12} holds true as well.
In fact, applying $\nabla_\xi$ to \eqref{xg_11}-\eqref{xg_21} we have the following two main possibilities:

\noindent
1) $\nabla_\xi$ can hit the profile $\hat{f_{\leq s^{1/8}}}$, therefore causing a loss of $s^{1/8}$; or

\noindent
2) $\nabla_\xi$ can hit the phase, which will introduce a growing factor of $s^{1+\delta_N}$.

\noindent
Since we have a bound of $s^{-7/4}$ on the $L^2$ norms of \eqref{xg_11}-\eqref{xg_21}, the desired bound of $s^{-1/2}$
for the $L^2$ norms of  $\nabla_\xi$  \eqref{xg_11}-\eqref{xg_21} follows. 

\subsection{Decay estimate for $e^{it\Lambda} G_2$} 
We write 
\begin{align*}
e^{it\Lambda} G_2 (t,x) & = 
\int_1^t e^{i(t-s)\Lambda} 
		\FF^{-1}_\xi \left[ \int_{\R^3} |\xi| e^{is \tilde{\psi} (\xi,\eta)} \hat{f_{\geq s^{1/8} }} (\xi-\eta,s) \overline{ \hat{f_{\geq s^{1/8} }} } (\eta,s) d\eta  \right] \, ds
\end{align*}
where $\tilde{\psi} (\xi,\eta) = {|\xi-\eta|}^2 - {|\eta|}^2 = {|\xi|}^2 - 2 \xi \cdot \eta$.
We start by using the symbol ${|\xi|}$ to integrate by parts in $\eta$.
Then by symmetry we can reduce to consider the following term:
\begin{equation}
\label{g_2decay}
\int_1^t e^{i(t-s)\Lambda}  \frac{1}{s} \, \FF^{-1}_\xi \left[ \int_{\R^3} \frac{\xi}{|\xi|} \,
			e^{is \tilde{\psi} (\xi,\eta)} \,  \nabla_\eta \hat{f_{\geq s^{1/8} }} (\xi-\eta,s) \overline{ \hat{f_{\geq s^{1/8} }} } (\eta,s) d\eta \right] 
			ds \, .
\end{equation}
The contribution of the time integral between $t-1$ and $t$ can be easily estimated by Sobolev embedding. To estimate the contribution of the integral from $1$ to $t-1$, we use the linear dispersive estimate for the wave equation and our large frequency cutoff convention to bound it as
(here $\vec R$ denotes the Riesz transform):
\begin{align*}
&\int_1^{t-1} \frac{1}{t-s} \,  \frac{1}{s} \, 
		\sum_{k = - \infty}^{ \log s^{\delta_N} } 2^{2k} {\left\|  P_k \vec R
		\left(  e^{is\Delta} x f_{\geq s^{1/8}}  \,  e^{is\Delta} f_{\geq s^{1/8}}  \right)  \right\|}_{L^1} \, ds
\\
 & \lesssim \int_1^{t-1} \frac{1}{t-s} \,  \frac{1}{s} \,  s^{2\delta_N}  
 						{\left\| e^{is\Delta} x f_{\geq s^{1/8}} \, e^{is\Delta} f_{\geq s^{1/8}} \right\|}_{L^1} \, ds
\\
 & \lesssim \int_1^{t-1} \frac{1}{t-s} \,  \frac{1}{s} \,  s^{2\delta_N} 
 						{\left\|  x f \right\|}_{L^2} {\| f_{\geq s^{1/8}}  \|}_{L^2} \, ds
\\
 & \lesssim \int_1^{t-1} \frac{1}{t-s} \,  \frac{1}{s} \,  s^{2\delta_N}  \, s^\d \frac{1}{s^{1/8}} s^\d  \, ds \lesssim \frac{1}{t} \, ,
\end{align*}
provided $2\delta_N + 2\d < \frac{1}{8}$. \endproof


\section{Improved $L^2$ estimates for $G$}\label{sec:imp}
In this section we show how to obtain some improved weighted $L^2$ estimates (Lemma \ref{Lemmag_1})
and some improved low frequency estimates (Lemma \ref{Lemmag_2}) for two different components of $G$. 
This will be essential in closing the a priori estimate on ${\| x^2 f \|}_{L^2}$ in section \ref{secweightedF}.
We split $G(f,\bar{f})$ into the two components
\begin{align*}
g_1 & :=   G(f_{\leq s^{1/4} } , \overline{f}_{\leq s^{1/4}} )
\\
g_2 & :=   G(f_{\geq s^{1/4} } , \overline{f} )  +   G(f_{\leq s^{1/4}} , \overline{f}_{\geq s^{1/4} } )  \, .
\end{align*}
The two terms in the definition of $g_2$ can be treated similarly, so we reduce to considering $g_1$ and $g_2$ given by
\begin{align}
\label{g_1}
\hat{g}_1 & := \int_0^t \int_{\R^3} |\xi| e^{is \psi (\xi,\eta)} 
			\hat{f_{\leq s^{1/4} }} (\xi-\eta,s) \overline{ \hat{f_{\leq s^{1/4} }} }(\eta,s) d\eta ds
\\
\label{g_2}
\hat{g}_2 & := \int_0^t \int_{\R^3} |\xi| e^{is \psi (\xi,\eta)} \hat{f_{\geq s^{1/4} }} (\xi-\eta,s) \overline{ \hat{f} } (\eta,s) d\eta ds \, .
\end{align}
Thanks to the small spatial support of the inputs in $g_1$ one can show that the ${\| x^2 g_1 \|}_{L^2}$ 
grows slower than ${\| x^2 g \|}_{L^2}$. In particular the following is true:
\begin{lemma}
 \label{Lemmag_1}
Let $g_1$ be defined by \eqref{g_1}, then
\begin{equation}
\label{x^2g_1}
{\| x^2 g_1 \|}_{L^2} \lesssim t^{3/4} \, .
\end{equation}
\end{lemma}

\proof
We apply as usual $\nabla_\xi^2$ to $\hat{g}_1$ and obtain:
\begin{subequations}
\begin{align}
\label{dxi^2g_11}
& \int_0^t \int_{\R^3} e^{i s \psi (\xi,\eta)} |\xi| \nabla_\xi^2 
			\hat{f_{\leq s^{1/4}} }(\xi-\eta,s) \overline{\hat{f_{\leq s^{1/4}}}} (\eta,s) d\eta ds
\\
\label{dxi^2g_12}
& \int_0^t \int_{\R^3} s \nabla_\xi \psi  e^{i s \psi (\xi,\eta)} |\xi|  \nabla_\xi 
			\hat{f_{\leq s^{1/4}}}(\xi-\eta,s) \overline{\hat{f_{\leq s^{1/4}}}} (\eta,s) d\eta ds
\\
\label{dxi^2g_13}
& \int_0^t \int_{\R^3} s^2 {(\nabla_\xi \psi )}^2 e^{i s \psi (\xi,\eta)} |\xi| 
			\hat{f_{\leq s^{1/4}} }(\xi-\eta,s) \overline{\hat{f_{\leq s^{1/4}} }} (\eta,s) d\eta ds
\end{align}
\end{subequations}
plus similar or easier terms. 
The first contribution is estimated by
\begin{align*}
{\| \eqref{dxi^2g_11} \|}_{L^2} & \lesssim \int_0^t s^{\delta_N} {\|  x^2 f_{\leq s^{1/4}} \|}_{L^2} {\|  e^{is\Delta} f_{\leq s^{1/4}}  \|}_{L^\infty} \, ds
\\
& \lesssim \int_0^t s^{\delta_N}  s^{1/4} {\|  x f_{\leq s^{1/4}} \|}_{L^2}  \frac{1}{s^{3/2}} {\|  f_{\leq s^{1/4}} \|}_{L^1} \, ds 
\\
& \lesssim \int_0^t s^{\delta_N}  s^{1/4} s^\d  \frac{1}{s^{3/2}} s^{1/8} s^\d  \, ds \lesssim 1 \, .
\end{align*}
The second term \eqref{dxi^2g_12} is bounded by
\begin{align*}
{\| \eqref{dxi^2g_12} \|}_{L^2} & \lesssim \int_0^t s^{1+2\delta_N} {\|  x f_{\leq s^{1/4}} \|}_{L^2} {\|  e^{is\Delta} f_{\leq s^{1/4}}  \|}_{L^\infty} \, ds
\\
& \lesssim \int_0^t s^{1+2\delta_N}  s^\d  \frac{1}{s^{3/2}} s^{1/8} s^\d  \, ds \lesssim t^{5/8 + 2(\d + \d_N)}
\end{align*}
which is okay for $\d$ and $\d_N$ small enough.
To estimate \eqref{dxi^2g_13} we integrate by parts in $\eta$ using \eqref{symg}.
This produces two terms of the same type as \eqref{dxi^2g_12} and can therefore be estimated similarly.
The validity of \eqref{x^2g_1} follows $_\Box$

\begin{remark} Notice that without the information on the small support of the inputs of $g_1$,
one could obtain, essentially just by using \eqref{symg} and H\"{o}lder, a bound of the form:
\begin{equation}
\label{x^2G}
{\| x^2 G \|}_{L^2} \lesssim t^{1 - \alpha + C\d} \, 
\end{equation}
which is much worse than \eqref{x^2g_1}.
Such a bound fails to suffice when trying to estimate $x^2 F$.
In particular, a term like \eqref{B_11},
with $G$ in place of $g_1$, would not be bounded a priori by $t^{1-2\a -\d}$.
\end{remark}

The following lemma gives an improved small frequency bound on $g_2$.
\begin{lemma}
\label{Lemmag_2}
Let $g_2$ be the term defined in \eqref{g_2}, and $P_{\leq k}$ denote the usual projection on frequencies smaller than $2^k$. 
Then
\begin{equation}
\label{P_kg_2}
{\| P_{\leq k} g_2 \|}_{L^2} \lesssim 2^{(7/4-3\d) k} \, .
\end{equation}
\end{lemma}

\proof
We start by using Bernstein's inequality twice, then H\"{o}lder,
and eventually the support information on $f$ to improve the decay in $s$ of the integrand:
\begin{align*}
& {\left\| P_{\leq k}  \int_0^t \Lambda  e^{-is |\nabla| } \left( e^{is\Delta} f_{\geq s^{1/4} }   e^{-is\Delta} \overline{ f } \right) \, ds \right\|}_{L^2}
\\
& \lesssim 2^k  \int_0^t {\left\| P_{\leq k} \left( e^{is\Delta} f_{\geq s^{1/4} }   e^{-is\Delta} \overline{ f } \right) \right\|}_{L^2} \, ds 
\\
& \lesssim 2^{ k } 2^{ (3/4-3\d)k }  \int_0^t {\left\|  e^{is\Delta} f_{\geq s^{1/4} }   e^{-is\Delta} \overline{ f } \right\|}_{ L^{{(3/4 - \d)}^{-1} }} \, ds 
\\
& \lesssim  2^{ (7/4-3\d) k}  \int_0^t {\|  f_{\geq s^{1/4} } \|}_{L^2}   {\| e^{-is\Delta} \overline{ f } \|}_{ L^{ {(1/4 - \d)}^{-1} } } \, ds 
\\
& \lesssim  2^{ (7/4-3\d) k }  \int_0^t s^\d \frac{1}{ {\langle s \rangle}^{1/4} }  \frac{1}{ {\l s \r}^{3/4 + 3\d}}  s^\d \, ds  \lesssim 2^{ (7/4-3\d)k }  \qquad   _\Box
\end{align*}

\vskip15pt
\section{Weighted estimates for the \S component} 
\label{secweightedF}

The purpose of this section is to prove: 
\begin{proposition}
\label{proweightedF}
Let $F$ be the bilinear term defined in \eqref{F}, then 
\begin{equation*}
t^{-\d} {\| x F \|}_{L^2} + t^{-1+2\a+\d} {\|  x^2 F \|}_{L^2}  \lesssim {\| (u,w) \|}_X^2 \, .
\end{equation*}
\end{proposition}
\noindent
A key identity that we are going to use is
\begin{equation}
\label{dxiphi0}
\nabla_\xi \phi = - 2\eta = - 2 \frac{\eta}{|\eta|}(\frac{\eta}{|\eta|}\cdot \nabla_\eta \phi)-2\frac{\phi}{|\eta|}\frac{\eta}{|\eta|}
\, .
\end{equation}
This is saying that $\nabla_\xi \phi$ vanishes (up to some mild singularity) on the resonant set 
$\mathscr{R_{\phi}} = \{ \phi = 0, \nabla_\eta \phi=0 \}$.
Therefore \eqref{dxiphi0} can be considered as a type of null structure, see again \cite{PS}.
In particular, we can use the factor $\nabla_\eta \phi$ to integrate by parts in $\eta$, and the factor $\phi$ to integrate by parts in $s$.


\begin{remark}[Simplification of \eqref{dxiphi0}]
\label{remdxiphi}
We can disregard the factors of $\eta/|\eta|$ which multiply $\nabla_\eta \phi$ and $\phi/|\eta|$ in \eqref{dxiphi0}.
Indeed, their presence is inconsequential whenever one is estimating $L^p$ norms of $e^{is\Lambda} g$ for $ 1 < p < \infty$. In the cases when $p = \infty$ is needed, we invoke the bound on the $\dot{B}^0_{\infty,1}$ norm of $e^{is\Lambda} g$.
Also, whenever an integration by parts in $\eta$ is performed using $\nabla_\eta \phi$,
one would need to consider the case when $\nabla_\eta$ hits the factor $\eta/|\eta|$.
Via Hardy's inequality all such terms are analogous to the terms obtained when $\nabla_\eta$ hits $\hat{g}$.
Therefore, to simplify the presentation, we will abuse notation and rephrase identity \eqref{dxiphi0} into:
\begin{equation}
\label{dxiphi}
\nabla_\xi \phi = \nabla_\eta \phi + \frac{\phi}{|\eta|} \, .
\end{equation}
\end{remark}

\subsection{Estimate ${\| xF \|}_{L^2} \lesssim t^\d$}
Applying $\nabla_\xi$ to the bilinear term $\hat{F}$ we get the following terms:
\begin{subequations}
\begin{align}
\label{dxif1}
& \int_0^t \int_{\R^3} e^{is \phi (\xi,\eta)} \nabla_\xi \hat{f}(\xi-\eta,s) \hat g (\eta,s) d\eta ds
\\
\label{dxif2}
& \int_0^t \int_{\R^3} s \nabla_\xi \phi \, e^{is \phi (\xi,\eta)} \hat{f} (\xi-\eta,s) \hat g (\eta,s) d\eta ds \, .
\end{align}
\end{subequations}
Using \eqref{dxiphi} to integrate by parts in $\eta$ and $s$ in the term \eqref{dxif2}, we get the contributions:
\begin{subequations}
\begin{align}
\label{dxif21}
& \int_0^t \int_{\R^3} e^{is \phi (\xi,\eta)} \nabla_\eta \hat{f}(\xi-\eta,s) \hat{g} (\eta,s) d\eta ds
\\
\label{dxif22}
& \int_0^t \int_{\R^3} e^{is \phi (\xi,\eta)} \hat{f} (\xi-\eta,s)  \nabla_\eta \hat{g} (\eta,s) d\eta ds
\\
\label{dxif23}
& \int_{\R^3} t \, e^{it \phi (\xi,\eta)} \hat{f} (\xi-\eta,s) \frac{1}{|\eta|} \hat{g} (\eta,s) d\eta
\\
\label{dxif24}
& \int_0^t \int_{\R^3} s \, e^{is \phi (\xi,\eta)} \hat{f} (\xi-\eta,s)  \frac{1}{{|\eta|}} \partial_s  \hat{g} (\eta,s) d\eta
\\
\label{dxif25}
& \int_0^t \int_{\R^3} s \, e^{is \phi (\xi,\eta)}  \partial_s \hat{f} (\xi-\eta,s)  \frac{1}{{|\eta|}} \hat{g} (\eta,s) d\eta \, .
\end{align}
\end{subequations}


\subsubsection{Estimate of \eqref{dxif1} and \eqref{dxif21}}
These terms are identical and can be bounded by a simple $L^2 \times L^\infty$ estimate\footnote{
Notice that here we are implicitely using the stronger bound ${\| e^{it\Lambda} g \|}_{\dot{B}^0_{\infty,1}} \lesssim t^{-1}$,
since factors of $\eta/|\eta|$ should appear in \eqref{dxif21}.
}:
\begin{align*}
{\| \eqref{dxif1} \|}_{L^2}
		& \lesssim  \int_0^t {\| x f \|}_{L^2} {\| e^{is\Lambda} g \|}_{L^\infty}  \, ds
			\lesssim \int_0^t s^\d \frac{1}{s} \, ds \lesssim t^\d \, .
\end{align*}

\subsubsection{Estimate of \eqref{dxif22}}
Recall that 
\begin{equation}
\label{gn_1}
g(t,x) = g (0,x) + G(t,x) = n_0(x) + i \Lambda^{-1} n_1 (x) + G(t,x) \, ,
\end{equation}
and therefore
\begin{align}
\label{xgn_1}
x g(t,x)  =  x  n_0(x) + i \Lambda^{-1} x n_1(x) - i \Lambda^{-2} \vec{R} n_1(x)  +  x G(x,t)  \, .
\end{align}
The contribution in \eqref{dxif22} coming from $ x n_0 + \Lambda^{-1} x n_1 + x G$
can be estimated easily by an $L^\infty \times L^2$ estimate using \eqref{decayu},
Hardy's inequality in combination with \eqref{data}, and the a priori bound on $x G$ in $L^2$.
To estimate the remaining contribution 
\begin{align}
\label{dxif22a}
 &  \int_0^t \int_{\R^3} e^{is \phi (\xi,\eta)} \hat{f} (\xi-\eta,s)  \frac{1}{{|\eta|}^2} \hat{ \vec{R}  n_1} (\eta,s) d\eta ds
\end{align}
we use instead an $L^2 \times L^\infty$ estimate, 
and the linear dispersive estimate \eqref{linearwave0}:
\begin{align*}
{\| \eqref{dxif22a} \|}_{L^2}
		& \lesssim  \int_0^t  {\| f \|}_{L^2} {\| e^{is\Lambda} \Lambda^{-2} \vec{R} n_1 \|}_{L^\infty}  
			\lesssim  \int_0^t  \, \frac{1}{ \l s \r}  {\| n_1 \|}_{\dot{B}^0_{1,1}}  \lesssim \log \l t \r \, .
\end{align*}

\subsubsection{Estimate of \eqref{dxif23}}
According to \eqref{gn_1} we distiguish the contibution coming from $n_0 + G$ and the one coming from $\Lambda^{-1} n_1$.
Using Hardy's inequality, and again \eqref{linearwave0}, we see that
\begin{align*}
{\| \eqref{dxif23} \|}_{L^2}  
		& \lesssim  t {\| e^{it\Delta} f \|}_{L^\infty} {\| e^{it\Lambda}  \Lambda^{-1} ( n_0 + G) \|}_{L^2}  
						 +  t {\| f \|}_{L^2} {\| e^{it\Lambda} \Lambda^{-2} \vec{R} n_1 \|}_{L^\infty} 
\\
& \lesssim t \, \frac{1}{t^{1+\a}}\| x (n_0  + G) \|_{L^2}
+ t \, \frac{1}{t} {\| n_1 \|}_{\dot{B}^0_{1,1}} \lesssim 1 \, .
\end{align*}

\subsubsection{Estimate of \eqref{dxif24}}
Since $e^{is\Lambda} \partial_s g = \Lambda {|u|}^2$, we have
${\| \Lambda^{-1}  e^{is\Lambda} \partial_s g \|}_{L^3} \lesssim {\| u \|}_{L^6}^2 \lesssim s^{-2+2\d}$.
Thanks to this and an $L^6 \times L^3$ estimate one sees that \eqref{dxif24} is bounded by $1$.

\subsubsection{Estimate of \eqref{dxif25}}
Finally, using $e^{is\Lambda} \partial_s f = u w$
we can bound
\begin{align*}
{\| \eqref{dxif25} \|}_{L^2}
		& \lesssim  \int_0^t s {\| e^{is\Delta} \partial_s f \|}_{L^\infty} {\| e^{is\Lambda} \Lambda^{-1} (n_0 + G) \|}_{L^2} 
						+ s {\| \partial_s f \|}_{L^2} {\| e^{is\Lambda} \Lambda^{-2} \vec{R} n_1 \|}_{L^\infty}  \, ds
\\
		& 	\lesssim \int_0^t s  {\| u \|}_{L^\infty} {\| w \|}_{L^\infty} \| x (n_0  + G) \|_{L^2}
				+ s {\| u \|}_{L^\infty} {\| w \|}_{L^2} \frac{1}{s} {\| n_1 \|}_{\dot{B}^0_{1,1}}  \, ds 
\\
		& 	\lesssim \int_0^t s \frac{1}{ {\l s \r}^{2+\a}} 
				\, ds \lesssim 1 \, .
\end{align*}

\subsection{Estimate ${\| x^2 F \|}_{L^2} \lesssim t^{1 - 2\a - \d}$}
Applying $\nabla_\xi^2$ to $\hat{F}$ gives the following three types of contributions
\begin{subequations}
\begin{align}
\label{dxi^2f1}
& \int_0^t \int_{\R^3} e^{is \phi(\xi,\eta)} \nabla_\xi^2 \hat{f}(\xi-\eta,s) \hat g (\eta,s) d\eta ds
\\
\label{dxi^2f2}
& \int_0^t \int_{\R^3} s \, \eta  \, e^{is \phi (\xi,\eta)} \nabla_\xi \hat{f} (\xi-\eta,s) \hat g (\eta,s) d\eta ds
\\
\label{dxi^2f3}
& \int_0^t \int_{\R^3} s^2 \eta^2  e^{is \phi (\xi,\eta)} \hat f(\xi-\eta,s) \hat g (\eta,s) d\eta ds \, .
\end{align}
\end{subequations}

\subsubsection{Estimate of \eqref{dxi^2f1}}
The term \eqref{dxi^2f1} can be easily estimated by an $L^2 \times L^\infty$ application of H\"{o}lder's inequality, so we skip it.

\subsubsection{Estimate of \eqref{dxi^2f2}}
For the term \eqref{dxi^2f2} we use \eqref{dxiphi} to integrate by parts in $\eta$ and $s$.
Notice that here one does not get a term containing $x^2 g$, for which we do not have good enough control.
Using \eqref{dxiphi} one gets the following types of contributions:
\begin{subequations}
\begin{align}
\label{dxi^2f21}
& \int_0^t \int_{\R^3} e^{is \phi(\xi,\eta)} \nabla^2_\eta \hat{f}(\xi-\eta,s) \hat{g} (\eta,s) d\eta ds
\\
\label{dxi^2f22}
& \int_0^t \int_{\R^3} e^{is \phi(\xi,\eta)} \nabla_\eta \hat{f} (\xi-\eta,s)  \nabla_\eta \hat{g} (\eta,s) d\eta ds
\\
\label{dxi^2f23}
& \int_{\R^3} t \frac{1}{|\eta|} \, e^{it \phi (\xi,\eta)} \nabla_\eta \hat{f} (\xi-\eta,s) \hat{g} (\eta,s) d\eta
\\
\label{dxi^2f24}
& \int_0^t \int_{\R^3} s \frac{1}{|\eta|} \, e^{is \phi (\xi,\eta)} \nabla_\eta \hat{f} (\xi-\eta,s)  \partial_s  \hat{g} (\eta,s) d\eta
\\
\label{dxi^2f25}
& \int_0^t \int_{\R^3} s \frac{1}{|\eta|} \, e^{is \phi (\xi,\eta)}  \partial_s \nabla_\eta \hat{f} (\xi-\eta,s)  \hat{g} (\eta,s) d\eta
\end{align}
\end{subequations}
plus simlar or easier terms. 

\eqref{dxi^2f21} is identical to \eqref{dxi^2f1} and can be estimated in the same way.

To estimate \eqref{dxi^2f22} we use \eqref{xgn_1} and distinguish again two different cases according to \eqref{xgn_1}.
Using H\"{o}lder's inequality, \eqref{SL^6}, and Sobolev's embedding,
the contribution coming from $x (n_0 + G) + i\Lambda^{-1} x n_1$ can be bounded by
\begin{align*}
		&  \int_0^t {\| e^{is\Delta} x f \|}_{L^6} {\left\| e^{is\Lambda} [ x (n_0 + G)  +  i\Lambda^{-1} x n_1 ] \right\|}_{L^3}  \, ds
		\\
		\lesssim & \int_0^t \frac{1}{s} {\| x^2 f \|}_{L^2} \left( {\| x (n_0 + G) \|}_{\dot{H}^\frac{1}{2}}
						 +  {\| \Lambda^{-{1/2}} x n_1 \|}_{L^2}  \right) \, ds
		\\
		\lesssim & \int_0^t \frac{1}{s} s^{1-2\a-\d} \, ds \lesssim t^{1-2\a-\d} \, .
\end{align*}
The contribution coming from the term which contains $e^{is\Lambda} \Lambda^{-2} \vec{R} n_1$ can be bounded
via an $L^2 \times L^\infty$ estimate similar to the one performed on the term \eqref{dxif22a}, so we skip it.

The term \eqref{dxi^2f23} can be treated similarly to \eqref{dxi^2f22}, 
since $\hat{g} / |\eta|$ plays the same role as $\nabla_\eta \hat{g}$,
and the factor of $t$ plays the same role of the integral in time.

\eqref{dxi^2f24} can be bounded by an $L^2 \times L^\infty$ estimate using
\begin{equation*}
{\left\| e^{is\Lambda} \Lambda^{-1} \partial_s g \right\|}_{L^\infty} = {\| u^2 \|}_{L^\infty} \lesssim \frac{1}{ {\l s \r}^{2+2\a}} \, .
\end{equation*}

The last term \eqref{dxi^2f25} is more delicate. To estimate it we need the following inequalities
\begin{align}
\label{dxi^2f25a}
& {\left\| e^{is\Delta} \partial_s x f \right\|}_{L^p} \lesssim \frac{1}{s^{4/3 - 2/p}} s^{\d + \delta_N} \qquad p = 2, 6 \, ,
\\
\label{dxi^2f25b}
& {\left\| e^{is\Lambda} \Lambda^{-1} (n_0 + G) \right\|}_{L^3} \lesssim \frac{1}{s^{1/3}} s^{\d +\delta_N} \, .
\end{align}
Postponing for the moment the proof of these, we bound \eqref{dxi^2f25} as follows:
\begin{align*}
{\| \eqref{dxi^2f25} \|}_{L^2}
		& \lesssim  \int_0^t s {\| e^{is\Delta} \partial_s x f \|}_{L^6} {\| e^{is\Lambda} \Lambda^{-1} (n_0 + G) \|}_{L^3}  
					+  s {\| e^{is\Delta} \partial_s x f \|}_{L^2} {\| e^{is\Lambda} \Lambda^{-2} \vec{R} n_1  \|}_{L^\infty}  \, ds
		\\
		& \lesssim  \int_0^t s \frac{1}{s} s^{\d + \delta_N} \frac{1}{s^{1/3}} s^{\d + \delta_N }  
					 + s \frac{1}{s^{1/3}} s^{\d + \delta_N} \frac{1}{s} {\| n_1 \|}_{\dot{B}^0_{1,1}} \, ds  
			\lesssim t^{2/3 + 2(\delta_N + \d)}
\end{align*}
which is majorized by $ t^{1-2\a-\d} = t^{2/3 + 3\d}$ provided $2\delta_N \leq \d$.

To conclude the bound on \eqref{dxi^2f2} we need to show \eqref{dxi^2f25a} and \eqref{dxi^2f25b}.
To see why \eqref{dxi^2f25a} holds observe that
\begin{align}
\label{dsxf}
e^{is\Delta} \partial_s x f & = \FF^{-1} \left( \int_{\R^3}  s \nabla_\xi \phi \, \hat{u}(\xi-\eta) \hat{w}(\eta) d\eta \right)
+  w \,  e^{is\Delta} x f \, .
\end{align}
In the case $p=2$ we can bound the first summand with an $L^6 \times L^3$ estimate to obtain the desired bound of $s^{\d + \d_N - 1/3}$.
For $p=6$ we can use instead an $L^6 \times L^\infty$ estimate to obtain a bound of $s^{\d + \d_N - 1}$.
The second summand in \eqref{dsxf} is easier to treat, so we skip it.
%
%
The bound \eqref{dxi^2f25b} is verified for the initial data $n_0$ 
since the linear dispersive estimate \eqref{linearwave} gives
\begin{align*}
& {\left\| e^{is\Lambda} \Lambda^{-1} n_0 \right\|}_{L^3} \lesssim \frac{1}{s^{1/3}} {\| \Lambda^{-1/3} n_0 \|}_{ L^{3/2} } 
		\lesssim  \frac{1}{s^{1/3}} {\| n_0 \|}_{L^{9/7}} 
		\lesssim  \frac{1}{s^{1/3}}  {\| \langle x \rangle  n_0 \|}_{L^2} \, ,
\end{align*}
having used Hardy-Littlewood-Sobolev
for the second inequality.
Moreover, we see that
\begin{equation*}
e^{is\Lambda} \Lambda^{-1} G (s) = \int_0^s e^{i(s-r) \Lambda} {| u (r) |}^2 \, dr \, ,
\end{equation*}
hence
\begin{align*}
{\left\| e^{is\Lambda} \Lambda^{-1} G \right\|}_{L^3} & \lesssim  \int_0^s \frac{1}{ {(s-r)}^{1/3} } {\| \Lambda^{2/3} u^2 \|}_{ L^{3/2} } \, dr
		\lesssim \int_0^s \frac{1}{ {(s-r)}^{1/3} } {\| u \|}_{ W^{1,3} } {\| u \|}_{L^3} \, dr
		\\
		& \lesssim \int_0^s \frac{1}{ {(s-r)}^{1/3}} \, r^{\delta_N} \frac{1}{\sqrt{ \l r \r }} \, r^{\d/2}  \frac{1}{\sqrt{ \l r \r}} \, r^{\d/2} \, dr
		\lesssim \frac{1}{s^{1/3}} s^{\d +\delta_N} \, .
\end{align*}

\vskip10pt
\subsubsection{Estimate of \eqref{dxi^2f3}}
Let us denote \eqref{dxi^2f3} by
\begin{equation}
\label{dxi^2f}
B(f,g) (t,\xi) = \int_0^t \int_{\R^3} s^2 \eta^2  e^{is \phi (\xi,\eta)} \hat f(\xi-\eta,s) \hat g (\eta,s) d\eta ds \, .
\end{equation}
To estimate this term one would be tempted to do the integration by parts algebra using \eqref{dxiphi}. Note however that this would lead to a term
containing $x^2 g$, and the available bound \eqref{x^2G} on this latter would not allow us to close the desired a priori estimate, 
depsite the more than integrable decay of $u$.
We then split $g$ as 
$ g = g_0 + g_1 + g_2 $,
where $g_0$ is the initial data and $g_1$ and $g_2$ are as in \eqref{g_1} and \eqref{g_2} respectively.
We split accordingly  $B(f,g)$ into $B_1$ and $B_2$ with
\begin{subequations}
\begin{align}
B_1 (f,g) &:= B(f,g_0 + g_1)
\\
B_2 (f,g) &:= B(f,g_2) \, .
\end{align}
\end{subequations}

\vskip10pt
\paragraph{\bf{Estimate of $B_1$ in $L^2$}}
We begin by looking at the contribution containing the initial data $g_0$, that is $B(f,g_0)$.
From \eqref{dxiphi0} and \eqref{dxiphi} we have  $\eta^2 \sim \phi + \eta \nabla_\eta \phi$.
Using this idenity to integrate by parts in time and frequency, one sees that $B(f,g_0)$ is given by the following main terms
\begin{subequations}
\begin{align}
\label{B_101}
& \int_0^t \int_{\R^3} s \eta e^{is \phi (\xi,\eta)} \nabla_\eta \hat{f} (\xi-\eta,s) \hat{g}_0 (\eta) d\eta ds
\\
& \label{B_102}
\int_0^t \int_{\R^3}  s \eta e^{is \phi (\xi,\eta)} \hat{f} (\xi-\eta,s) \nabla_\eta \hat{g}_0 (\eta) d\eta ds
\\
& \label{B_103}
\int_{\R^3} t^2  e^{it \phi (\xi,\eta)} \hat{f} (\xi-\eta,t) \hat{g}_0 (\eta) d\eta \, 
\\
& \label{B_104}
\int_0^t \int_{\R^3} s^2  e^{it \phi (\xi,\eta)} \partial_s \hat{f} (\xi-\eta,s) \hat{g}_0 (\eta) d\eta \,ds .
\end{align}
\end{subequations}
The first term \eqref{B_101} is analogous to \eqref{dxi^2f2} and can be treated in the same way, so we can skip it.
To estimate \eqref{B_102} we notice that 
$\eta \nabla_\eta \hat{g}_0 = \eta \nabla_\eta \hat{n}_0 + \nabla_\eta \hat{n}_1 + \hat{n}_1 / {|\eta|}$,
One can then use the linear dispersive estimate \eqref{linearwave} and obtain
\begin{align*}
{\| e^{is\Lambda} \FF^{-1} (\eta \nabla_\eta \hat{g}_0 )\|}_{L^3} \lesssim \frac{1}{s^{1/3}} s^{\d_N}\, .
\end{align*}
The above estimate and an $L^6 \times L^3$ application of H\"{o}lder's inequality
show the desired bound 
${\| \eqref{B_102} \|}_{L^2} \lesssim t^{2/3 + \d + \d_N}$.
\eqref{B_103} is easily estimated by
\begin{align*}
{\| \eqref{B_103} \|}_{L^2} & \lesssim  t^2  {\|  e^{is\Delta} f \|}_{L^6}  {\|  e^{it\Lambda} g_0  \|}_{L^3}
\lesssim  t^2 \frac{1}{t} t^\d  \frac{1}{ t^{1/3} }  
\lesssim  t^{2/3  + \d} \lesssim t^{1-2\a -\d}\, ,
\end{align*}
where we used \eqref{aprioriG} and \eqref{param}.
Using once again $e^{is\Delta} \partial_s f = u w$ we obtain 
\begin{align*}
{\| \eqref{B_104} \|}_{L^2} & \lesssim \int_0^t s^2 {\|  e^{is\Delta} \partial_s f \|}_{L^6}  {\|  e^{is\Lambda} g_0  \|}_{L^3} \, ds
\\
& \lesssim \int_0^t s^2 {\| u \|}_{L^6} {\| w \|}_{L^\infty} \frac{1}{s^{1/3}} \, ds 
\lesssim \int_0^t s^2 \frac{1}{s} s^\d \frac{1}{s} \frac{1}{s^{1/3}} \, ds \lesssim  t^{2/3  + \d} \lesssim t^{1-2\a -\d} \, 
\end{align*}
by \eqref{aprioriG} and $\eqref{param}$.

We are now left with the contribution in $B_1$ coming from
\begin{equation}
\label{B_1}
B (f,g_1) =  \int_0^t \int_{\R^3} s^2 \eta^2  e^{is \phi (\xi,\eta)} \hat{f} (\xi-\eta,s) \hat{g}_1 (\eta,s) d\eta ds \, .
\end{equation}
In order to prove that this term satisfies the desired bound of $t^{1 - 2\a - \d}$ 
we proceed again by using \eqref{dxiphi} to integrate by parts in $\eta$ and $s$.
The terms obtained by doing so are of the type \eqref{dxi^2f21}--\eqref{dxi^2f25} (or easier),
or they are the analogue of \eqref{B_103} and \eqref{B_104} with $g_1$ instead of $g_0$,
except for the following two terms:
\begin{subequations}
\begin{align}
\label{B_11}
& \int_0^t \int_{\R^3} e^{is \phi (\xi,\eta)} \hat{f} (\xi-\eta,s) \nabla_\eta^2 \hat{g}_1 (\eta,s) d\eta ds \, ,
\\
& \label{B_12}
\int_0^t \int_{\R^3} s^2 \phi \, e^{is \phi (\xi,\eta)} \hat{f} (\xi-\eta,s) \partial_s \hat{g}_1 (\eta,s) d\eta ds \, .
\end{align}
\end{subequations}
%
%
Thanks to Lemma \ref{Lemmag_1} we have
\begin{align*}
{\| \eqref{B_11} \|}_{L^2} & \lesssim \int_0^t {\|  e^{is\Delta} f \|}_{L^\infty}  {\|  x^2 g_1  \|}_{L^2} \, ds
\\
& \lesssim \int_0^t \frac{1}{s^{1+\a}} s^{3/4} \, ds  \lesssim  t^{3/4 -\a} \, .
\end{align*}
This is majorized by $t^{1-2\a-\d}$ provided $\a + \d \leq \frac{1}{4}$, which is consistent with the choice \eqref{param}.
Finally, \eqref{B_12} can be bounded in a straightforward fashion by an $L^6 \times L^3$ estimate,
using ${\| e^{is\Lambda} \partial_s g_1 \|}_{L^3} = {\| \Lambda (e^{is\Delta} f_{\leq s^{1/4}})^2 \|}_{L^3} \lesssim s^{\delta_N} \|e^{is\Delta}f_{\leq s^{1/4}}\|_{L^6}^{2} \lesssim s^{-2 + 2\d + \d_N}$.


\vskip10pt
\paragraph{\bf{Estimate of $B_2$ in $L^2$}}
To estimate $B_2$ we decompose it further according to the size of the frequency $\eta$.
Let $\chi$ be a smooth positive radial and compactly supported function which equals $1$ on $[0,1]$ and vanishes on $[2,\infty)$,
and define $\chi_{\leq K}= \chi(\frac{\cdot}{K})$ and $\chi_{\geq K} = 1 - \chi_{\leq K}$.
Let $l$ be a positive number to be determined later, define 
\begin{align}
\label{dxi^2fl}
B_2^{\mbox{\tiny low} } (f,g) (t,\xi) & := 
				\int_0^t \int_{\R^3} s^2 \eta^2 \chi_{ \leq s^{-l}} (\eta) e^{is \phi (\xi,\eta)} \hat f(\xi-\eta,s) \hat g_2(\eta,s) d\eta ds
\\
\label{dxi^2fh}
B_2^{\mbox{\tiny high}} (f,g) (t,\xi) & := \int_0^t \int_{\R^3} s^2 \eta^2 \chi_{\geq s^{-l}} (\eta)
	 e^{is \phi (\xi,\eta)} \hat f(\xi-\eta,s) \hat g_2(\eta,s) d\eta ds \, .
\end{align}


\vskip5pt
\subparagraph{{\it Estimate of $B_2^{\mbox{\tiny low} }$}}
The term $B_2^{\mbox{\tiny low} }$ can be treated directly by using the smallness of the symbol and Lemma \ref{Lemmag_2}.
Applying H\"{o}lder's and and Bernstein's inequalities, we obtain
\begin{align*}
{\| B_2^{\mbox{\tiny low} } (f,g) \|}_{L^2}
		& \lesssim \int_0^t s^2 \frac{1}{s^{2l}}  {\| e^{is\Delta} f \|}_{L^6} {\| P_{ \leq \log_2 (s^{-l}) } e^{-is \Lambda} g_2 \|}_{L^3} \, ds 
\\
& \lesssim \int_0^t s^2 \frac{1}{ s^{2l} } \frac{1}{s} s^\d  \frac{1}{ s^{\frac{l}{2}} } {\| P_{ \leq \log_2 (s^{-l}) } g_2 \|}_{L^2} \, ds 
\\
& \lesssim \int_0^t s \frac{1}{ s^{(\frac{17}{4} - 3\d)l} }  s^\d \, ds \lesssim  t^2 \frac{1}{ t^{ \frac{17}{4} l} } t^{(1+3l)\d} \, .
\end{align*}
We then choose 
\begin{equation}
\label{cut}
l = \frac{1}{3} - \frac{1}{60}
\end{equation}
in such a way that the resulting bound at the end of the above chain of inequalities is majorized by $t^{2/3 + 2\d}$,
which is less than $t^{1-2\a-\d}$ as desired.

\vskip5pt
\subparagraph{{\it Estimate of $B_2^{\mbox{\tiny high} }$}}
To estimate the component $B_2^{\mbox{\tiny high} }$ in \eqref{dxi^2fh} we use once more \eqref{dxiphi} to integrate by parts in time and frequency.
By doing this one obtains again terms of the form \eqref{dxi^2f21}--\eqref{dxi^2f25}  (or easier ones),
or the analogues of \eqref{B_103}--\eqref{B_104} with $g_2$ instead of $g_0$, 
plus the following term:
\begin{equation}
\label{dxi^2fh1}
\int_0^t \int_{\R^3} \chi_{\geq s^{-l}} (\eta) e^{is \phi_{\pm}(\xi,\eta)} \hat{f}(\xi-\eta,s) \nabla_\eta^2 \hat{g_2}(\eta,s) d\eta ds \, .
\end{equation}
Notice once again that here we do not have access to a good estimate on $x^2 g_2$.
However we can use the fact that $|\eta|$ is not too small,
and use the available bound on the $L^2$ norm $\Lambda x^2 g$.
Using H\"{o}lder's inequality, followed by Sobolev's embedding and Bernstein's inequality we can estimate
\begin{align*}
{\| \eqref{dxi^2fh1} \|}_{L^2}
		& \lesssim \int_0^t {\| e^{is\Delta} f \|}_{L^6} {\| P_{ \geq \log_2 (s^{-l}) } e^{-is \Lambda} \left( x^2 g_2 \right) \|}_{L^3} \, ds 
\\
& \lesssim \int_0^t  \frac{1}{s} s^\d   {\| P_{ \geq \log_2 (s^{-l}) } \Lambda^\frac{1}{2} x^2 g_2 \|}_{L^2} \, ds 
\\
& \lesssim \int_0^t \frac{1}{s}  s^\d s^{l/2} {\| \Lambda x^2 g_2 \|}_{L^2} \, ds  \lesssim  t^{l/2+\d} t^{1-3\a}  \, ,
\end{align*}
where we used the fact that $\|\Lambda x^2 g_2\|_{L^2}$ 
satisfies the same bounds as $\|\Lambda x^2 g\|_{L^2}$,
because applying the spatial localizations on $f$ in the definition of $g_2$ does not affect the arguments in Section \ref{lamdax2G}.
This gives the desired bound of $t^{1-2\a-\d}$ provided we can choose $\a$ and $l$ such that
\begin{equation*}
\a \geq \frac{l}{2} + 2\d \, .
\end{equation*}
For  $\a$ given by \eqref{param} and $l$ given by \eqref{cut}
this inequality holds true provided $\d \leq \frac{1}{480}$.
This concludes the proof of a priori estimates on $x^2 F$ and hence of Proposition \ref{proweightedF}.
Together with Propositions \ref{proenergy}, \ref{proweightedG} and \ref{prodecayw}, this gives
the desired apriori bound for solutions of \eqref{Z} as explained in Section \ref{secnorms}, from which 
Theorem \ref{maintheo} follows.  \endproof

\vskip10pt
\subsection*{Acknowledgements}
The first two authors were supported in part by the Simons Postdoctoral Fellows Program. The authors would like to thank the referee for his careful reading of the manuscript and all the helpful comments and suggestions.

\end{document}